\documentclass[onefignum,onetabnum]{siamart171218}

\usepackage{lipsum}
\usepackage{amsfonts}
\usepackage{graphicx}
\usepackage{epstopdf}
\usepackage{algorithmic}
\ifpdf
  \DeclareGraphicsExtensions{.eps,.pdf,.png,.jpg}
\else
  \DeclareGraphicsExtensions{.eps}
\fi

\newsiamremark{remark}{Remark}
\newsiamremark{hypothesis}{Hypothesis}
\crefname{hypothesis}{Hypothesis}{Hypotheses}
\newsiamthm{claim}{Claim}

\headers{An algebraic geometry perspective on TDA}{Paul Breiding}

\title{An algebraic geometry perspective on \\ topological data analysis}

\author{Paul Breiding\thanks{Technische Universit\"at Berlin
  (\email{p.breiding@tu-berlin.de},
  \url{page.math.tu-berlin.de/\textasciitilde breiding/}),
  \funding{The author has received funding from the European Research Council under the European Union’s Horizon 2020 research and innovation programme (grant agreement No 787840)}.}
}

\usepackage{amsopn}

\begin{document}

\maketitle

The story of Topological Data Analysis (TDA) is undoubtedly a story of success that includes a wide range of diverse applications. To contribute to this story and to shed light on it from a new perspective, I will survey TDA from the point of view of \emph{Algebraic Geometry}.

First of all, there are some immediate applications of algebra and algebraic geometry to TDA.
The data structure behind \emph{Persistent Homology} (PH) is the persistence module, which is an inherently algebraic concept. Furthermore, attempts to extend PH to multiple parameters use concepts from Commutative Algebra and Algebraic Combinatorics \cite{CZ2009, HOST2017, Miller2017}. Here, however, I would like to discuss some less well-known roles that algebraic geometry can play for TDA. Specifically, I would like to present applications of Numerical Algebraic Geometry (NAG) and Enumerative Algebraic Geometry (EAG).

NAG is concerned with computing numerical solutions to a system of $n$ polynomial equations $F(x)=(f_1(x),\ldots, f_n(x))=0$ in $n$ variables $x=(x_1,\ldots,x_n)$ over the complex numbers.
The textbook \cite{SW} is a standard reference. The computational paradigm in NAG is \emph{numerical homotopy continuation}. The idea behind this is to generate a system of equations $G(x)$ of which the solutions are known (a so called \emph{start-system}) and then continue the solutions of $G(x)=0$ along a deformation of $G(x)$ towards $F(x)$.
The continuation leads to an ordinary differential equation, called a \emph{Davidenko-differential equation}, which is solved by standard numerical predictor-corrector methods for ODEs. The state of the art implementations of homotopy continuation are \texttt{Bertini}~\cite{BHSW06}, \texttt{HOM4PS}~\cite{HOM4PS}, \texttt{HomotopyContinuation.jl}~\cite{BT}, \texttt{NAG4M2}~\cite{Leykin2018} and \texttt{PHCPack}~\cite{PHCPack}.

EAG, on the other hand, counts the number of solutions of a system of polynomial equations, often by the means of intersection theory. The textbook \cite{EH} provides a good introduction to those fields. Although different at first sight, NAG and EAG are intimately related: The key benefit of NAG is that one can generate initial values \emph{for all isolated solutions} of $F(x)=0$ in $\mathbb{C}^n$. For instance, if the degree of the $i$-th polynomial is $d_i$, then $G(x)=(x_i^{d_i}-a_i)_{i=1}^n$, where $a_1,\ldots,a_n\in\mathbb C^*$, may serve as a start system for the homotopy $(1-t)G(x)+tF(x)$, $0\leq t\leq 1$. The number of solutions of $G(x)=0$ is $D=d_1\cdots d_n$ and a theorem from algebraic geometry implies that the number of isolated solutions of $F(x)=0$ is \emph{at most} $D$.
Therefore, continuing the solutions of $G(x)=0$ towards $F(x)=0$ produces all isolated solutions of $F(x)=0$. In practice, however, the number of solutions of $F(x)=0$ is significantly smaller than $D$ and diverging solutions must be eliminated. EAG helps to construct other start systems, adapted to the structure of $F(x)$, that improve the efficiency of the algorithm. The article \cite{BST2020} explains this relation between NAG and EAG in detail.

Back to TDA: Let us consider the situation in which the model $M\subset \mathbb{R}^n$ is given as the zero set of $s$ polynomials in $n$ variables $F(x)=(f_1(x),\ldots, f_s(x))$. In algebraic geometry such an $M$ is called a \emph{real algebraic variety}. As an example, consider the conformation space of cyclooctane. Cyclooctane is a molecule built out of eight carbon atoms $x_1,\ldots,x_8\in \mathbb{R}^3$ that are aligned in a ring such that the distances between neighboring atoms are equal to a constant $c>0$. The energy of a configuration $(x_1,\ldots,x_8)$ is minimized when the angles between successive bonds are all equal to $\arccos(-\frac{1}{3}) \approx 109.5^\circ$. The polynomial equations in $3\cdot 8 = 24$
variables are
\begin{align*}&\Vert x_1 - x_2\Vert^2 = \cdots = \Vert x_7 - x_8\Vert^2 = \Vert x_8 - x_1\Vert^2  = c^2\\
&\Vert x_1 - x_3 \Vert^2 =\cdots = \Vert x_6 - x_8 \Vert^2 = \Vert x_7 - x_1 \Vert^2 =\Vert x_8 - x_2 \Vert^2 = \frac{8}{3}c^2.
\end{align*}
It is known that the solution set of these equations, up to simultaneous translation and rotation, is homeomorphic to a union of the Klein bottle and a sphere, which intersect in two rings \cite{CMTW2010}. How can one use TDA for computing this?

\begin{figure}[h]
\begin{center}
\includegraphics[height = 5cm]{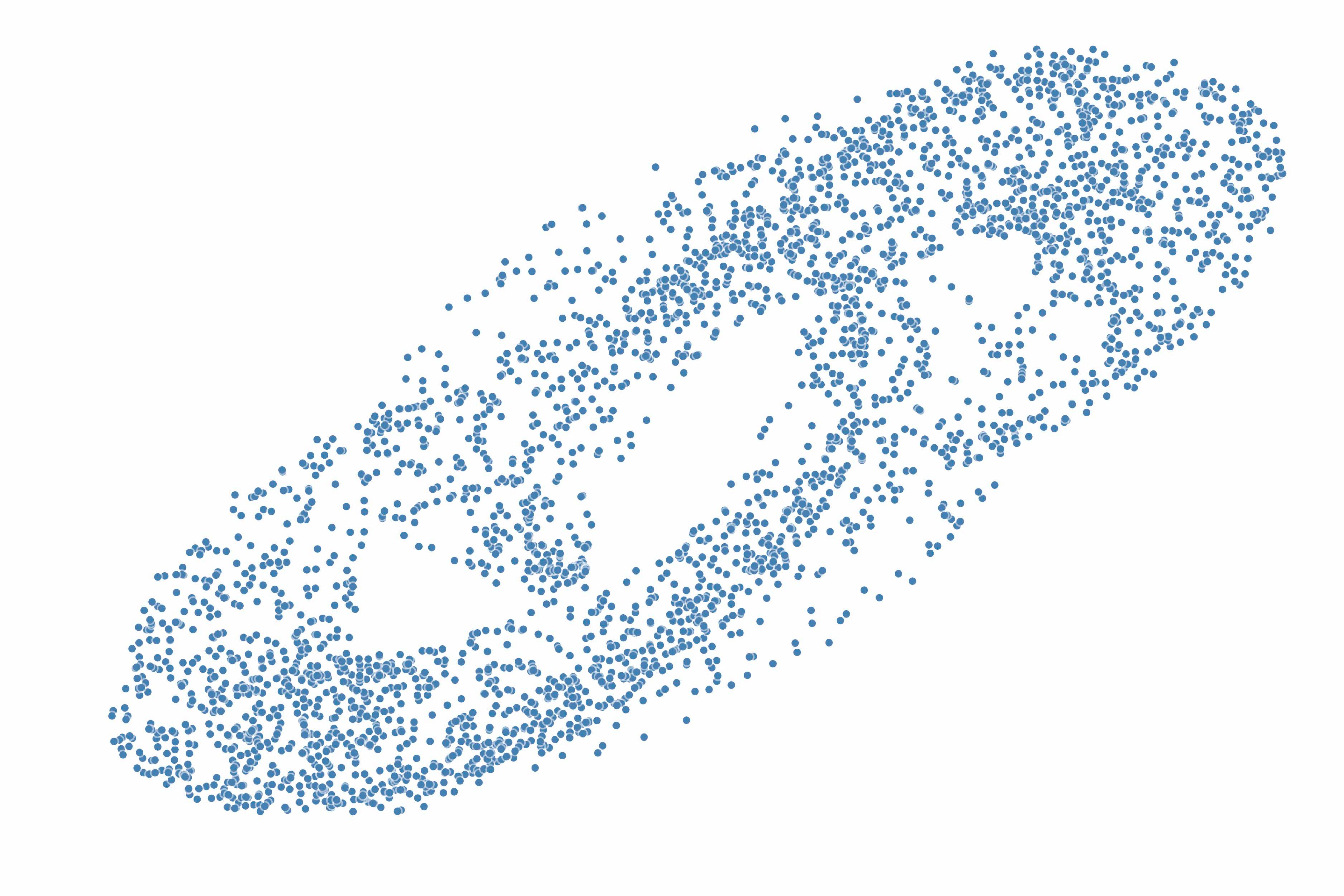}
\end{center}
\caption{\label{fig1} The picture shows a sample from the cyclooctane variety for $c^2=2$ projected to a two-dimensional space.Due to translational and rotational invariance $x_1=(0,0,0)$ and $x_2=(c,0,0)$ are fixed and the last entry of $x_3$ is set equal to zero. See \cite{cyclooctane2019} for the code that produced this data.
}
\end{figure}

NAG can be used for generating sample of points from $M$, which may then serve as input for persistent homology; figure \ref{fig1} shows a sample from the cyclooctane variety. One idea is to sample linear spaces $L$ of dimension equal to the codimension of $M$ and compute the points in $M\cap L$. Another idea is to sample points $q\in\mathbb{R}^n$ in the ambient space and compute the point on $M$ that minimizes the distance to~$q$. Both computational problems can be cast as a system of polynomial equations and solved using NAG. Research directions in NAG include how to sample with respect to a probability distribution on $M$ \cite{BO2019} and produce samples that are as dense in $M$ as desired \cite{DEHH2018}.

NAG and EGA are also used to study two important numbers for TDA: the homological feature size $\mathrm{hfs}(M)$ \cite{CEH2007} and the reach $\tau(M)$ \cite{NSW2008} of $M$. Horobet and Weinstein \cite{HW2018} showed that, if $M$ is an algebraic manifold (i.e., a real algebraic variety which is also a manifold) defined by polynomials over $\mathbb Q$, then both $\mathrm{hfs}(M)$ and $\tau(M)$ are algebraic over $\mathbb Q$. Therefore, both can be computed by the means of NAG. I discuss this for the reach.

The reach of $M$ is $\tau(M)=\min\{\tfrac{1}{\sigma(M)}, \tfrac{1}{2}\rho(M)\}$, where $\sigma(M)$ is the maximal curvature of a geodesic running through $M$ and $\rho(M)$ is the width of the narrowest \emph{bottleneck} of $M$ ($\tau(M)$ can also be defined as the distance from $M$ to its medial axis; see, e.g., \cite{ACKMRW2019}, for an explanation of the equivalence of the definitions). A bottleneck is a pair $(x,y)\in M^2$, $x\neq y$, such that $x-y$ is perpendicular to both tangent spaces~$\mathrm{T}_x M$ and~$\mathrm{T}_yM$. Provided $\mathrm{codim}(M) = s$, this can be cast as a system of polynomial equations in $2n+2s$ variables:
$$B(x,y,\lambda,\mu) = \begin{bmatrix} F(x) &  F(y) & JF(x) \lambda - (x-y) & JF(y) \mu - (x-y)\end{bmatrix} = 0,$$
where $\lambda = (\lambda_1,\ldots,\lambda_s)$, $\mu = (\mu_1,\ldots,\mu_s)$ and $JF$ is the Jacobian matrix of $F$. This system of equations can be solved using NAG and the real solutions can be extracted from the complex solutions. Recall that in NAG one computes the \emph{isolated} solutions of a polynomial system -- the trivial solutions for which $x=y$ are not computed! Bottlenecks have been studied intensely: Eklund \cite{Eklund2018} discussed bottlenecks from the perspective of NAG and di~Rocco et.\ al.\ \cite{DEW2019} gave a formula for the number of (complex) solutions of $B(x,y,\lambda,\mu)=0$ in terms of \emph{polar classes} of $M$.

Equations for $\sigma(M)$ are less straightforward. For planar curves, however, there is a direct formula for the curvature $\gamma(x)$ at $x\in M$ such that $\sigma(M)=\min_{x\in M}\gamma(x)$. In this case, the first-order optimality equations for $\sigma(M)$, that the gradient of $\gamma(x)$ is perpendicular to the tangent space $\mathrm{T}_x M$, give a system of polynomial equations. Solving this yields $\sigma(M)$.

\begin{figure}[h]
\begin{center}
\includegraphics[width = 0.475\textwidth]{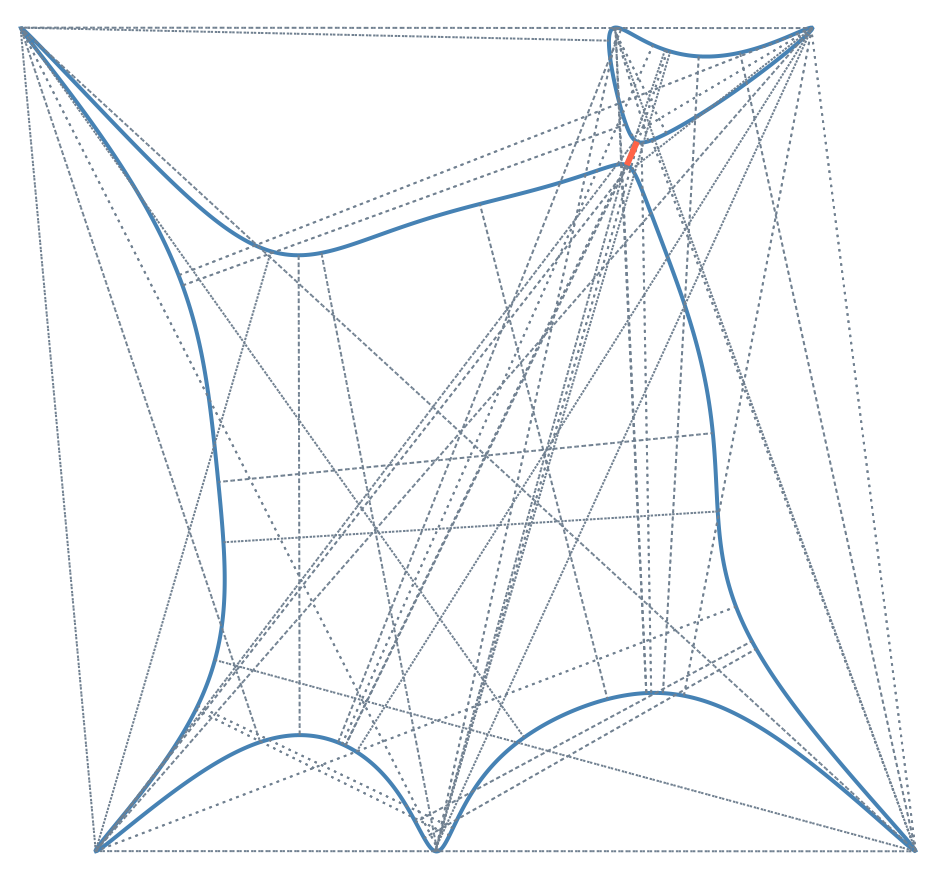}\includegraphics[width = 0.475\textwidth]{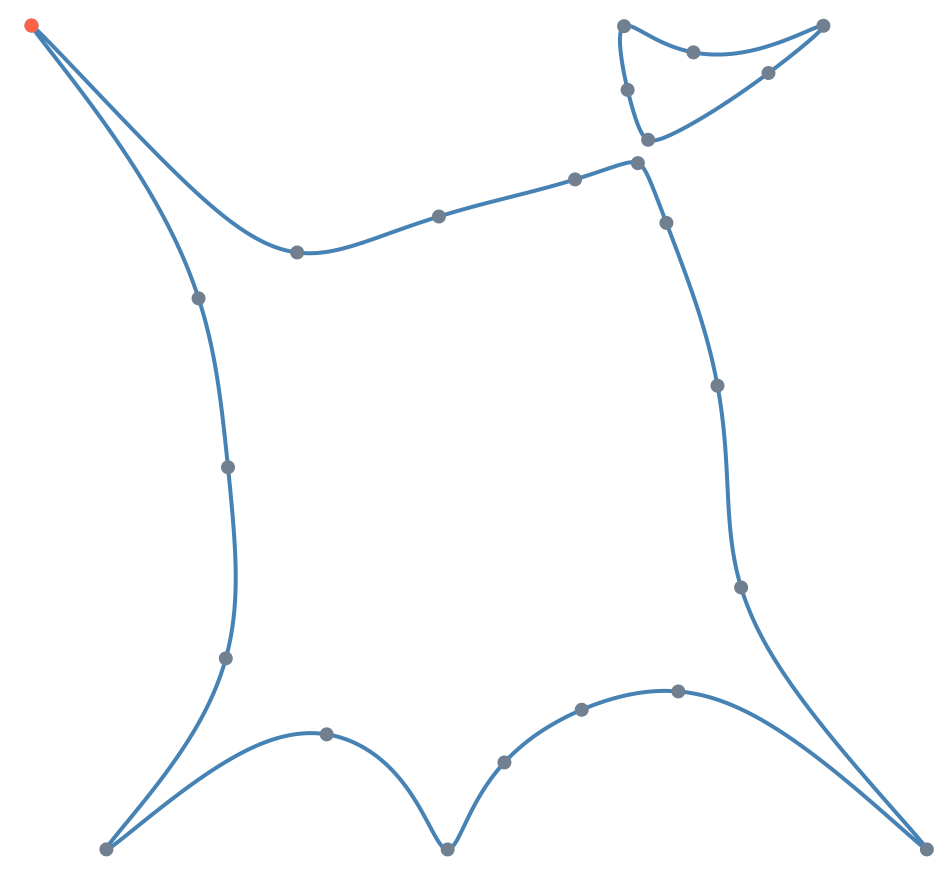}
\end{center}
\caption{\label{fig2} The pictures depict the ingredients for the computation of the reach of the planar curve $C=\{(x^3 - xy^2 + y + 1)^2(x^2 + y^2 - 1) + y^2 = 5\}$. The left picture displays all bottlenecks of~$C$. The narrowest bottleneck is red and its width is $\approx 0.138$. The right picture shows all points of critical curvature of $C$. The red point is the point of maximal curvature $\approx 2097.17$. Therefore, $\tau(C) \approx \{\tfrac{1}{2097.17}, \tfrac{0.138}{2}\} = \tfrac{1}{2097.17}$. The code that produced these results is available at \cite{reach2019}.}
\end{figure}

Summarizing, one can compute $\rho(M)$ and $\sigma(M)$ separately using NAG, and the reach $\tau(M)$ can be inferred from those numbers. Figure \ref{fig2} shows an example of this computation for a planar curve. It should be mentioned that, if $\tau(M)$ is known, \cite{NSW2008} provides a criterion for when the correct homology of $M$ can be computed from a finite sample of $M$. Combined with the sampling algorithms from above this yields provably correct algorithms for computing homology.

I also want to mention the research in \cite{BCT2019, BCL2019, CKS2018}. In these articles the reach is replaced by a lower bound that involves the \emph{real condition number} of a system of polynomials.
This lower bound holds not only for real algebraic varieties but also for the more general class of \emph{semialgebraic sets}. Using the condition number the authors derive a complexity analysis of an algorithm for computing homology.

In concluding this survey I want to propose three possible future directions in the field of algebraic geometry for TDA. The first is an analysis in the sense of NAG and EAG for $\sigma(M)$ -- like \cite{Eklund2018, DEW2019} did for bottlenecks. This is indispensable in computing the reach beyond planar curves. The second direction is considering persistent homology using ellipsoids instead of balls: The experiments in  \cite{BKSW2018} show that this approach can greatly improve the quality of the output diagrams in persistent homology, yet a theoretical explanation is missing. The third direction is about sampling. The standard approach to sampling from nonlinear objects is by using \emph{Markov Chain Monte Carlo} (MCMC) methods. Combining this approach with NAG seems promising.

\bibliographystyle{siamplain}

\end{document}